\documentclass{article}%
\usepackage{amssymb}
\usepackage{amsmath}%
\usepackage{amsfonts}%
\usepackage{graphicx}
\newtheorem{theorem}{Theorem}

\newtheorem{lemma}[theorem]{Lemma}

\newtheorem{proposition}[theorem]{Proposition}

\newenvironment{proof}[1][Proof]{\noindent\textbf{#1.} }{\ \rule{0.5em}{0.5em}}
\begin{document}

\title{Solitary waves in Abelian Gauge Theories}
\author{Vieri Benci$^{\ast}$, Donato Fortunato$^{\ast\ast}$\\$^{\ast}$Dipartimento di Matematica Applicata ''U. Dini''\\Universit\`{a} degli Studi di Pisa Universit\`{a} di Pisa\\Via Filippo Buonarroti 1/c, 56126 Pisa, Italy\\e-mail: benci@dma.unipi.it\\$^{\ast\ast}$Dipartimento di Matematica \\Universit\`{a} di Bari,\\Via Orabona 4, 70125 Bari, Italy\\e-mail: fortunat@dm.uniba.it}
\date{}
\maketitle

\begin{abstract}
Abelian gauge theories consist of a class of field equations which provide a
model for the interaction between matter and electromagnetic fields. In this
paper we analyze the existence of solitary waves for these theories. We assume
that the lower order term W is positive and we prove the existence of solitary
waves if the coupling between matter and electromagnetic field is small. We
point out that the positiveness assumption on W implies that the energy is
positive: this fact makes these theories more suitable to model physical phenomena.

\end{abstract}
\tableofcontents

\section{Introduction}

Roughly speaking a solitary wave is a solution of a field equation whose
energy travels as a localized packet and which preserves this localization in time.

In this paper we are interested in investigating the existence of solitary
waves relative to Abelian gauge theories. These theories consist of a class of
field equations that provide a model for the interaction of matter with the
electromagnetic field.

This problem has already been treated in \cite{bf}, \cite{tea}, \cite{tea2},
\cite{ca}, \cite{ingl}, \cite{cossu}, \cite{dav} and \cite{vaira}. In those
papers the lower order term $W$ has the form
\[
W(u)=\frac{1}{2}u^{2}-\frac{1}{p}\left|  u\right|  ^{p},\;\;2<p<6.
\]
In this paper, we analyze the case
\begin{equation}
W(u)\geq0:\label{due}%
\end{equation}
this assumption implies that the energy density is positive; the interest of
this feature will be discussed at the end of section \ref{uu} and in section
\ref{SW}.

\subsection{Abelian gauge theories\label{agt}}

Let $G$ be a subgroup of $U(N),$ the unitary group in $\mathbb{C}^{N},$ and
denote by $\Lambda^{k}(\mathbb{R}^{4},\mathfrak{g})$ the set of $k$-forms
defined in $\mathbb{R}^{4}$ with values in the Lie algebra $\mathfrak{g}$ of
the group $G$. A 1-form
\[
\Gamma=\sum_{j=0}^{3}\Gamma_{j}dx^{j}\in\Lambda^{1}(\mathbb{R}^{4}%
,\mathfrak{g})
\]
is called connection form. The operator
\[
d_{\Gamma}=\Lambda^{k}(\mathbb{R}^{4},\mathfrak{g})\rightarrow\Lambda
^{k+1}(\mathbb{R}^{4},\mathfrak{g})
\]
defined by
\[
d_{\Gamma}=d+\Gamma=\sum_{j=0}^{3}\left(  \frac{\partial}{\partial x^{j}%
}+\Gamma_{j}\right)  dx^{j}%
\]
is called covariant differential and the operators
\[
D_{j}=\frac{\partial}{\partial x^{j}}+\Gamma_{j}:\mathcal{C}^{1}\left(
\mathbb{R}^{4},\mathbb{C}^{N}\right)  \rightarrow\mathcal{C}^{0}\left(
\mathbb{R}^{4},\mathbb{C}^{N}\right)  ,\;\;\;j=0,...,3
\]
are called covariant derivatives. The 2-form
\[
F=d_{\Gamma}\Gamma=\sum_{i,j=0}^{3}\left(  \partial_{i}\Gamma_{j}+\left[
\Gamma_{i},\Gamma_{j}\right]  \right)  dx^{i}\wedge dx^{j}%
\]
is called \textit{curvature}.

Now, we equip $\mathbb{R}^{4}$ with the Minkowski quadratic form given by
\[
\left\langle v,v\right\rangle _{M}=-\left|  v_{0}\right|  ^{2}+\sum_{j=1}%
^{3}\left|  v_{j}\right|  ^{2}.
\]
where $v=(v_{0},v_{1},v_{2},v_{3})\;$is a 4-vector with components in
$\mathbb{R}$ or in $\mathbb{C.\;}$The Minkowski quadratic form can be extended
to the space of the differential forms $\alpha\in\Lambda^{k}(\mathbb{R}%
^{4},\mathfrak{g}),$ and it will be denoted by $\left\langle \alpha
,\alpha\right\rangle _{M}$ We set
\[
\mathcal{L}_{0}=-\frac{1}{2}\left\langle d_{\Gamma}\psi,d_{\Gamma}%
\psi\right\rangle _{M},
\]
where $\psi\in\mathbb{C}^{N}$.

$.$ We set
\[
\text{ }\mathcal{L}_{1}=-\frac{1}{2q^{2}}\left\langle d_{\Gamma}%
\Gamma,d_{\Gamma}\Gamma\right\rangle _{M},
\]
where $q>0$ is a real parameter which controls the coupling of $\mathcal{L}%
_{1}$ with $\mathcal{L}_{0}.$

A gauge field by definition (see e.g.\cite{yangL}, \cite{rub}), is a critical
point of the action functional
\begin{equation}
\mathcal{S}=\int\mathcal{L\;}dxdt,\;\;\mathcal{L}=\mathcal{L}_{0}%
+\mathcal{L}_{1}-W(\psi),\label{completa}%
\end{equation}
where $W:\mathbb{C}^{N}\rightarrow\mathbb{R}$ is a function which is assumed
to be $G$-invariant, namely
\begin{equation}
W(g\psi)=W(\psi),\ \;\;g\in G\mathbf{.}\label{4}%
\end{equation}

We are interested in the Abelian gauge theory, namely in the case in which
\[
G=U(1)=S^{1}=\left\{  z\in\mathbb{C:}\left|  z\right|  =1\right\}
\]
In this case the $\Gamma_{j}(x,t)$ are imaginary numbers and
\[
\left[  \Gamma_{i},\Gamma_{j}\right]  =0
\]
Then
\[
\mathcal{L}_{1}=-\frac{1}{2q^{2}}\left\langle d_{\Gamma}\Gamma,d_{\Gamma
}\Gamma\right\rangle _{M}=-\frac{1}{2q^{2}}\left\langle d\Gamma,d\Gamma
\right\rangle _{M}%
\]
If we set%

\[
A^{j}=A_{j}=-\frac{1}{iq}\Gamma_{j},\;\;j=1,2,3
\]
and
\[
\varphi=A^{0}=-A_{0}=\frac{1}{iq}\Gamma_{0},
\]

\noindent it turns out that $(A^{0},A^{1},A^{2},A^{3})\;$is a real valued
4-vector field and $\mathbf{A:}=(A^{1},A^{2},A^{3})$ is its spacial component.
Setting $x^{0}=t$, the covariant derivatives take the form
\[
D_{t}=\frac{\partial}{\partial t}+iq\varphi,\;\;\text{ }D_{j}=\frac{\partial
}{\partial x^{j}}-iqA_{j}\text{ }j=1,2,3\text{ }%
\]
and, for $q=0,$ they reduces to the usual ones. Using the above notation, the
Lagrangian density $\mathcal{L}_{0},$ can be written as follows
\begin{align}
\mathcal{L}_{0}  & =\frac{1}{2}\left|  D_{t}\psi\right|  ^{2}-\frac{1}%
{2}\left|  \mathbf{D}_{x}\psi\right|  ^{2}\\
& =\frac{1}{2}\left[  \left|  \left(  \frac{\partial}{\partial t}%
+iq\varphi\right)  \psi\right|  ^{2}-\left|  \left(  \nabla-iq\mathbf{A}%
\right)  \psi\right|  ^{2}\right]
\end{align}
where $\mathbf{D}_{x}\psi=\left(  D_{1}\psi,D_{2}\psi,D_{3}\psi\right)  $ and
$\mathcal{L}_{1}$ takes the form%

\[
\mathcal{L}_{1}=\frac{1}{2}\left|  \frac{\partial\mathbf{%
A%
}}{\partial t}+\nabla\varphi\right|  ^{2}-\frac{1}{2}\left|  \nabla
\times\mathbf{A}\right|  ^{2}.
\]
Here $\nabla\times$ and $\nabla$ denote respectively the curl and the gradient operators.

Making the variation of $\mathcal{S}$ with respect to $\psi,$ $\varphi$ and
$\mathbf{A}$ we get the following system of equations
\begin{equation}
D_{t}^{2}\psi-\mathbf{D}_{x}^{2}\psi+W^{\prime}(\psi)=0\label{e1+}%
\end{equation}
\begin{equation}
\nabla\cdot\left(  \frac{\partial\mathbf{%
A%
}}{\partial t}+\nabla\varphi\right)  =q\left(  \operatorname{Im}\frac{1}{\psi
}\frac{\partial\psi}{\partial t}+q\varphi\right)  \left|  \psi\right|
^{2}\label{e2+}%
\end{equation}
\begin{equation}
\nabla\times\left(  \nabla\times\mathbf{A}\right)  +\frac{\partial}{\partial
t}\left(  \frac{\partial\mathbf{%
A%
}}{\partial t}+\nabla\varphi\right)  =q\left(  \operatorname{Im}\frac
{\nabla\psi}{\psi}-q\mathbf{A}\right)  \left|  \psi\right|  ^{2}\label{e3+}%
\end{equation}
where
\[
\mathbf{D}_{x}^{2}\psi=\sum_{j=1}^{3}D_{j}^{2}\psi
\]
is the ''equivariant Laplacian''; also we have assumed the notation
\[
W^{\prime}(\psi)=\frac{\partial W}{\partial\psi_{1}}+i\frac{\partial
W}{\partial\psi_{2}},\text{ }\psi=\psi_{1}+i\psi_{2}%
\]

Notice that (\ref{4}) implies that $W$ is a functon of $\left|  \psi\right|  $
and we have
\begin{align*}
W(e^{i\vartheta}\psi)  & =W(\psi).\\
W^{\prime}(e^{i\vartheta}\psi)  & =e^{i\vartheta}W^{\prime}(\psi).
\end{align*}

In order to give a more meaningful form to the above equations, we will write
$\psi$ in polar form
\[
\psi(x,t)=u(x,t)\,e^{iS(x,t)},\;\;u\geq0,\;\;S\in\mathbb{R}/2\pi\mathbb{Z}%
\]
So (\ref{completa}) takes the following form
\begin{align*}
\mathcal{S(}u,S,\varphi,\mathbf{A}) &  =\int\int\left[  \frac{1}{2}\left(
\frac{\partial u}{\partial t}\right)  ^{2}-\frac{1}{2}\left|  \nabla u\right|
^{2}-W(u)\right]  dxdt+\\
&  +\frac{1}{2}\int\int\left[  \left(  \frac{\partial S}{\partial t}%
+q\varphi\right)  ^{2}-\left|  \nabla S-q\mathbf{A}\right|  ^{2}\right]
\,u^{2}dxdt\\
&  +\frac{1}{2}\int\int\left(  \left|  \frac{\partial\mathbf{A}}{\partial
t}\mathbf{+}\nabla\varphi\right|  ^{2}-\left|  \nabla\times\mathbf{A}\right|
^{2}\right)  dxdt
\end{align*}
and the equations (\ref{e1+},\ref{e2+},\ref{e3+}) take the form:
\begin{equation}
\square u+W^{\prime}(u)+\left[  \left|  \nabla S-q\mathbf{A}\right|
^{2}-\left(  \frac{\partial S}{\partial t}+q\varphi\right)  ^{2}\right]
\,u=0\label{e1}%
\end{equation}
\begin{equation}
\frac{\partial}{\partial t}\left[  \left(  \frac{\partial S}{\partial
t}+q\varphi\right)  u^{2}\right]  -\nabla\cdot\left[  \left(  \nabla
S-q\mathbf{A}\right)  u^{2}\right]  =0\label{e2}%
\end{equation}
\begin{equation}
\nabla\cdot\left(  \frac{\partial\mathbf{%
A%
}}{\partial t}+\nabla\varphi\right)  =q\left(  \frac{\partial S}{\partial
t}+q\varphi\right)  u^{2}\;\label{e3}%
\end{equation}
\begin{equation}
\nabla\times\left(  \nabla\times\mathbf{A}\right)  +\frac{\partial}{\partial
t}\left(  \frac{\partial\mathbf{%
A%
}}{\partial t}+\nabla\varphi\right)  =q\left(  \nabla S-q\mathbf{A}\right)
u^{2}\;\label{e4}%
\end{equation}

As we will see in the next section, these equations provide a model for Electrodynamics.

\subsection{The Maxwell equations\label{uu}}

In order to show the relation of the above equations with the Maxwell
equations and to get a model for Electrodynamics, we make the following change
of variables:
\begin{equation}
\mathbf{E=-}\left(  \frac{\partial\mathbf{%
A%
}}{\partial t}+\nabla\varphi\right) \label{pos1}%
\end{equation}
\begin{equation}
\mathbf{H}=\nabla\times\mathbf{A}\label{pos2}%
\end{equation}
\begin{equation}
\rho=-\left(  \frac{\partial S}{\partial t}+q\varphi\right)  qu^{2}\label{car}%
\end{equation}
\begin{equation}
\mathbf{j}=\left(  \nabla S-q\mathbf{A}\right)  qu^{2}\label{tar}%
\end{equation}
So (\ref{e3}) and (\ref{e4}) are the second couple of the Maxwell equations
with respect to a matter distribution whose charge and current density are
respectively $\rho$ and $\mathbf{j}$:
\begin{equation}
\nabla\cdot\mathbf{E}=\rho\label{gauss}%
\end{equation}
\begin{equation}
\nabla\times\mathbf{H}-\frac{\partial\mathbf{E}}{\partial t}=\mathbf{j}%
.\label{ampere}%
\end{equation}
(\ref{pos1}) and (\ref{pos2}) give rise to the first couple of the Maxwell
equation:
\begin{equation}
\nabla\times\mathbf{E+}\frac{\partial\mathbf{H}}{\partial t}=0\label{faraday}%
\end{equation}
\begin{equation}
\nabla\cdot\mathbf{H}=0.\label{monopole}%
\end{equation}
Equation (\ref{e1}) can be written as follows
\begin{equation}
\square u+W^{\prime}(u)+\frac{\mathbf{j}^{2}-\rho^{2}}{q^{2}u^{3}%
}=0\label{materia}%
\end{equation}
and finally Equation (\ref{e2}) is the charge continuity equation
\begin{equation}
\frac{\partial}{\partial t}\rho+\nabla\cdot\mathbf{j}=0.\label{continuità}%
\end{equation}

Notice that equation (\ref{continuità}) can be deduced by (\ref{gauss}) and
(\ref{ampere}). Thus, equations (\ref{gauss},..,\ref{materia}) are equivalent
to equations (\ref{e1},...,\ref{e4}). Concluding, an Abelian gauge theory, via
equations (\ref{gauss},..,\ref{materia}), provides a model of interaction of
the matter field $\psi$ with the electromagnetic field $(\mathbf{E}%
,\mathbf{H})$.

The gauge group is given by $\mathcal{G}\cong\mathcal{C}^{\infty}%
(\mathbb{R}^{4})$ with the additive structure and it acts on the variables
$\psi\mathbf{%
,%
}\varphi,\mathbf{%
A\;%
}$as follows
\begin{align*}
T_{\chi}\psi & =\psi e^{i\chi};\\
T_{\chi}\varphi & \rightarrow\varphi-\frac{\partial\chi}{\partial t};\\
T_{\chi}\mathbf{%
A%
}  & =\mathbf{%
A%
}+\nabla\chi;
\end{align*}
with $\chi\in\mathcal{C}^{\infty}(\mathbb{R}^{4}).\;$Equations (\ref{e1}%
,...,\ref{e4}) are gauge invariant by the way they have been
constructed.\ However, this fact, in eq. (\ref{gauss},..,\ref{materia}), can
be checked directly since the variables $u,\rho,\mathbf{j},\mathbf{E}%
,\mathbf{H}$ are themselves gauge invariant.

In order to show the relevance of assumption (\ref{due}), we will compute the
formula of the energy.

\begin{theorem}
If $(u,S,\mathbf{E},\mathbf{H})$ is a solution of the field equations
(\ref{e1},...\ref{e4}), its energy is given by
\[
\mathcal{E}(u,S,\mathbf{E},\mathbf{H})=\int\left[  \frac{1}{2}\left(
\frac{\partial u}{\partial t}\right)  ^{2}+\frac{1}{2}\left|  \nabla u\right|
^{2}+W(u)+\frac{\rho^{2}+\mathbf{j}^{2}}{2q^{2}u^{2}}\,+\frac{\mathbf{E}%
^{2}+\mathbf{H}^{2}}{2}\right]  dx
\]

\end{theorem}

\begin{proof}
We recall the well known expression for the energy density (see e.g.
\cite{Gelfand}):
\[
\frac{\partial\mathcal{L}}{\partial\left(  \frac{\partial u}{\partial
t}\right)  }\cdot\frac{\partial u}{\partial t}+\frac{\partial\mathcal{L}%
}{\partial\left(  \frac{\partial S}{\partial t}\right)  }\cdot\frac{\partial
S}{\partial t}+\frac{\partial\mathcal{L}}{\partial\left(  \frac{\partial
\varphi}{\partial t}\right)  }\cdot\frac{\partial\varphi}{\partial t}%
+\frac{\partial\mathcal{L}}{\partial\left(  \frac{\partial\mathbf{A}}{\partial
t}\right)  }\cdot\frac{\partial\mathbf{A}}{\partial t}-\mathcal{L.}%
\]

Now we will compute each term. We have:
\begin{equation}
\frac{\partial\mathcal{L}}{\partial\left(  \frac{\partial u}{\partial
t}\right)  }\cdot\frac{\partial u}{\partial t}=\left(  \frac{\partial
u}{\partial t}\right)  ^{2}\label{l1}%
\end{equation}
\begin{align*}
\frac{\partial\mathcal{L}}{\partial\left(  \frac{\partial S}{\partial
t}\right)  }\cdot\frac{\partial S}{\partial t}  & =\left(  \frac{\partial
S}{\partial t}+q\varphi\right)  \frac{\partial S}{\partial t}\,u^{2}\\
& =\left(  \frac{\partial S}{\partial t}+q\varphi\right)  \frac{\partial
S}{\partial t}\,u^{2}+\left(  \frac{\partial S}{\partial t}+q\varphi\right)
q\varphi u^{2}-\left(  \frac{\partial S}{\partial t}+q\varphi\right)  q\varphi
u^{2}\\
& =\left(  \frac{\partial S}{\partial t}+q\varphi\right)  ^{2}\,u^{2}-\left(
\frac{\partial S}{\partial t}+q\varphi\right)  q\varphi u^{2}\\
& =\frac{\rho^{2}}{q^{2}u^{2}}+\rho\varphi.
\end{align*}

In the following we write $\int g(x)dx$ instead of $\int_{\mathbb{R}^{3}%
}g(x)dx.$ By the Gauss equation (\ref{gauss}), multiplying by $\varphi$ and
integrating, we get
\[
-\int\mathbf{E\cdot}\nabla\varphi dx=\int\rho\varphi dx
\]

Thus, replacing this expression in the above formula, we get
\begin{equation}
\int\left(  \frac{\partial\mathcal{L}}{\partial\left(  \frac{\partial
S}{\partial t}\right)  }\cdot\frac{\partial S}{\partial t}\right)
dx=\int\left(  \frac{\rho^{2}}{q^{2}u^{2}}-\mathbf{E\cdot}\nabla
\varphi\right)  dx\label{l2}%
\end{equation}
Also we have
\begin{equation}
\frac{\partial\mathcal{L}}{\partial\left(  \frac{\partial\varphi}{\partial
t}\right)  }\cdot\frac{\partial\varphi}{\partial t}=0\label{l3}%
\end{equation}
and
\begin{equation}
\frac{\partial\mathcal{L}}{\partial\left(  \frac{\partial\mathbf{A}}{\partial
t}\right)  }\cdot\frac{\partial\mathbf{A}}{\partial t}=\left(  \frac
{\partial\mathbf{A}}{\partial t}\mathbf{+}\nabla\varphi\right)  \cdot
\frac{\partial\mathbf{A}}{\partial t}=-\mathbf{E}\cdot\frac{\partial
\mathbf{A}}{\partial t}\label{l4}%
\end{equation}
Moreover, using the notation (\ref{pos1},...,\ref{tar}), we have that
\[
\mathcal{L}=\frac{1}{2}\left(  \frac{\partial u}{\partial t}\right)
^{2}-\frac{1}{2}\left\vert \nabla u\right\vert ^{2}-W(u)+\frac{\rho
^{2}-\mathbf{j}^{2}}{2q^{2}u^{2}}\,+\frac{\mathbf{E}^{2}-\mathbf{H}^{2}}{2}%
\]

Then, by (\ref{l1},...,\ref{l4}) and the above expression for $\mathcal{L}$ we
get
\begin{align*}
\mathcal{E(}u,S,\varphi,\mathbf{A})  & =\int\left(  \frac{\partial\mathcal{L}%
}{\partial\left(  \frac{\partial u}{\partial t}\right)  }\cdot\frac{\partial
u}{\partial t}+\frac{\partial\mathcal{L}}{\partial\left(  \frac{\partial
S}{\partial t}\right)  }\cdot\frac{\partial S}{\partial t}+\frac
{\partial\mathcal{L}}{\partial\left(  \frac{\partial\mathbf{A}}{\partial
t}\right)  }\cdot\frac{\partial\mathbf{A}}{\partial t}-\mathcal{L}\right)
dx\\
& =\int\left(  \left(  \frac{\partial u}{\partial t}\right)  ^{2}+\frac
{\rho^{2}}{q^{2}u^{2}}-\mathbf{E\cdot}\nabla\varphi\mathbf{-E}\cdot
\frac{\partial\mathbf{A}}{\partial t}-\mathcal{L}\right)  dx\\
& =\int\left(  \left(  \frac{\partial u}{\partial t}\right)  ^{2}+\frac
{\rho^{2}}{q^{2}u^{2}}+\mathbf{E}^{2}\right)  dx\\
& -\int\left(  \frac{1}{2}\left(  \frac{\partial u}{\partial t}\right)
^{2}-\frac{1}{2}\left\vert \nabla u\right\vert ^{2}-W(u)+\frac{\rho
^{2}-\mathbf{j}^{2}}{2q^{2}u^{2}}\,+\frac{\mathbf{E}^{2}-\mathbf{H}^{2}}%
{2}\right)  dx\\
& =\int\left(  \frac{1}{2}\left(  \frac{\partial u}{\partial t}\right)
^{2}+\frac{1}{2}\left\vert \nabla u\right\vert ^{2}+W(u)+\frac{\rho
^{2}+\mathbf{j}^{2}}{2q^{2}u^{2}}\,+\frac{\mathbf{E}^{2}+\mathbf{H}^{2}}%
{2}\right)  dx
\end{align*}

\end{proof}

By the above theorem, the condition (\ref{due}) implies that the energy
density is positive and this fact makes these equations more suitable to model
physical phenomena.

\subsection{Solitary waves and solitons in the nonlinear wave
equation\label{SW}}

By \textit{solitary wave} we mean a solution of a field equation whose energy
travels as a localized packet; by \textit{soliton}, we mean a solitary wave
which exhibits some strong form of stability. This is a rather weak definition
of soliton but probably it is the most commonly used. Solitons have a
particle-like behavior and they occur in many questions of mathematical
physics, such as classical and quantum field theory, non linear optics, fluid
mechanics, plasma physics (see e.g. \cite{dodd}, \cite{fel}, \cite{Rajaraman},
\cite{Witham}).

In order to prove the existence of solitons, first it is necessary to prove
the existence of solitary waves and then to prove their stability.

We start our discussion of solitary waves considering eq. (\ref{e1}%
,..,\ref{e4}) in the case in which the covariant derivative is replaced by the
usual one. In this case the matter decouples with the Maxwell equations,
equation (\ref{e1+}) becomes the nonlinear wave equation
\begin{equation}
\frac{\partial^{2}\psi}{\partial t^{2}}-\triangle\psi+W^{\prime}%
(\psi)=0\label{KG}%
\end{equation}
and the Maxwell equations become homogeneous.

The easiest way to produce solitary waves of (\ref{KG}) consists in solving
the static equation
\begin{equation}
-\triangle u+W^{\prime}(u)=0\label{KGS}%
\end{equation}
and setting
\begin{equation}
\psi_{v}(t,x)=\psi_{v}(t,x_{1},x_{2},x_{3})=u\left(  \frac{x_{1}-vt}%
{\sqrt{1-v^{2}}},x_{2},x_{3}\right)  ;\label{pina}%
\end{equation}
$\psi_{v}(t,x)$ is a solution of eq. (\ref{KG}) which represents a bump which
travels in the $x_{1}$-direction with speed $v.$

In \cite{Poho} and \cite{strauss}, it has been proved that eq. (\ref{KGS}) has
nontrivial solutions provided that $W$ has the following form:
\begin{equation}
W(u)=\frac{1}{2}m^{2}u^{2}-\frac{1}{p}u^{p},\;m>0,\;2<p<6,\label{W1}%
\end{equation}

Moreover Shatah \cite{shatah} found a condition which guarantees the "orbital
stability" of the solitary waves of eq. (\ref{KG}); if $W$ is given by
(\ref{W1}), this condition becomes $2<p<\frac{10}{3}$ (see e.g. \cite{bergamo}
or \cite{sammomme}).

However, it would be interesting to assume $W\geq0;$ in fact the energy of a
solution of equation (\ref{KG}) is given by
\[
E(\psi)=\int\left[  \frac{1}{2}\left(  \frac{\partial\psi}{\partial t}\right)
^{2}+\frac{1}{2}\left|  \nabla\psi\right|  ^{2}+W(\psi)\right]  dx
\]
In this case, the positivity of the energy, not only is an important request
for the physical models related to this equation, but it provides good
\textit{a priori} estimates for the solutions of the relative Cauchy problem.
These estimates allow to prove the existence and well-posedness results under
very general assumptions on $W$.

Unfortunately Derrick \cite{derrick}, in a very well known paper, has proved
that request (\ref{due}) implies that equation (\ref{KGS}) has only the
trivial solution. His proof is based on the following equality (which in a
different form was found also by Pohozaev \cite{Poho}; for details see also
\cite{sammomme}). The Derrick-Pohozaev identity states that for any finite
energy solution $u$ of eq.(\ref{KGS}) it holds
\begin{equation}
\frac{1}{6}\int\left|  \nabla u\right|  ^{2}dx+\int W(u)dx=0\label{ventuno}%
\end{equation}
Clearly (\ref{ventuno}) and (\ref{due}) imply that $u\equiv0$.

However, we can try to prove the existence of solitons of eq. (\ref{KG}) (with
assumption (\ref{due})) exploiting the possible existence of \textit{standing
waves}, since this fact is not prevented by equation (\ref{ventuno}). A
\textit{standing wave} is a finite energy solution of (\ref{KG}) having the
following form
\begin{equation}
\psi_{0}(t,x)=u(x)e^{-i\omega_{0}t}\text{, }u\text{ }\geq0,\;\text{real}%
\label{sw}%
\end{equation}

\noindent Substituting (\ref{sw}) in eq.(\ref{KG}), we get
\begin{equation}
-\Delta u+W^{\prime}(u)=\omega_{0}^{2}u\label{static}%
\end{equation}

\noindent The Lagrangian of eq. (\ref{KG}) is given by
\[
\mathcal{L}=\frac{1}{2}\left(  \frac{\partial\psi}{\partial t}\right)
^{2}-\frac{1}{2}\left|  \nabla\psi\right|  ^{2}-W(\psi)
\]
It is easy to check that this Lagrangian is invariant for the Lorentz group.
Thus given a solution $\psi(t,x)$ of (\ref{KG}), we can obtain an other
solution $\psi_{1}(t,x)$ just making a Lorentz transformation on it. Namely,
if we take the velocity $\mathbf{v}=(v,0,0),$ $\left|  v\right|  <1$, and set
\[
t^{\prime}=\gamma\left(  t-vx_{1}\right)  ,\text{ }x_{1}^{\prime}%
=\gamma\left(  x_{1}-vt\right)  ,\text{ }x_{2}^{\prime}=x_{2},\text{ }%
x_{3}^{\prime}=x_{3}\;\;\;\text{with}\;\;\;\gamma=\frac{1}{\sqrt{1-v^{2}}}%
\]
it turns out that
\[
\psi_{\mathbf{v}}(t,x)=\psi(t^{\prime},x^{\prime})
\]
is a solution of (\ref{KG}).

In particular given a standing wave $\psi(t,x)=u(x)e^{-i\omega_{0}t},$ the
function $\psi_{\mathbf{v}}(t,x):=\psi(t^{\prime},x^{\prime})$ is a solitary
wave which travels with velocity $\mathbf{v.}$ Thus, if $u(x)=u(x_{1}%
,x_{2},x_{3})$ is any solution of Eq. (\ref{static}), then
\begin{equation}
\psi_{\mathbf{v}}(t,x_{1},x_{2},x_{3})=u\left(  \gamma\left(  x_{1}-vt\right)
,x_{2},x_{3}\right)  e^{i(\mathbf{k\cdot x}-\omega t)},\;\text{ }%
\label{solitone}%
\end{equation}
is a solution of Eq. (\ref{KG}) provided that%

\begin{equation}
\omega=\gamma\omega_{0}\;\;\text{and\ \ }\;\mathbf{k}=\gamma\omega
_{0}\mathbf{v}\label{33}%
\end{equation}
Notice that (\ref{pina}) is a particular case of (\ref{solitone}) when
$\omega_{0}=0.$

We have the following result:

\begin{theorem}
\label{lillo}Assume that

\begin{itemize}
\item (i) $W(u)\geq0,$

\item (ii) $W(0)=W^{\prime}(0)=0;\;W^{\prime\prime}(0)=m_{0}^{2}>0,$

\item (iii) there exists $u_{0}\in\mathbb{R}^{+},\;$ $W(u_{0})<\frac{1}%
{2}m_{0}^{2}u_{0}^{2}.$
\end{itemize}

\noindent Then eq. (\ref{KG}) has finite energy solitary waves of the form
$\;\psi_{0}(t,x)=u(x)e^{-i\omega_{0}t}$ for every frequency $\omega_{0}%
\in\left(  m_{1},m_{0}\right)  $ where
\[
m_{1}=\inf\left\{  m>0\mathbb{:\exists}u\in\mathbb{R}^{+},\;W(u)-\frac{1}%
{2}m^{2}u^{2}<0\right\}
\]

\end{theorem}

Notice that by (iii), $m_{1}<m_{0},$ then the interval $\left(  m_{1}%
,m_{0}\right)  $ is not empty.

\begin{proof}
By the previous discussion, it is sufficient to show that equation
(\ref{static}) has a solution $u$ with finite energy. The solutions of finite
energy of (\ref{static}) are the critical points in the Sobolev space
$H^{1}\left(  \mathbb{R}^{3}\right)  $ of the \textit{reduced action}
functional:
\begin{equation}
J(u)=\frac{1}{2}\int\left|  \nabla u\right|  ^{2}dx+\int G(u)dx,\text{
}G\left(  u\right)  =W(u)-\frac{1}{2}\omega_{0}^{2}u^{2}\text{ }%
\end{equation}

By a theorem of Berestycki and Lions \cite{BL81}, the existence of nontrivial
critical points of $J$ is guaranteed by the following assumptions on $G$:

\begin{itemize}
\item $G(0)=G^{\prime}(0)=0$

\item $G^{\prime\prime}(0)>0$

\item $\underset{s\rightarrow\infty}{\lim\sup}\frac{G^{\prime}(s)}{s^{5}}%
\geq0$

\item $\exists u_{0}\in\mathbb{R}^{+}:\;G(u_{0})<0.$
\end{itemize}

It is easy to check that for every frequency $\omega_{0}\in\left(  m_{1}%
,m_{0}\right)  ,$ the above assumptions are satisfied.
\end{proof}

\subsection{Solitary waves in Abelian gauge theories}

Now let us consider the problem of the existence of solitary waves for an
Abelian gauge theory. The Lagrangian $\mathcal{L}$ is invariant for the
following representation of the Lorentz group:
\begin{align*}
\psi_{\mathbf{v}}(t,x)  & =\psi(t^{\prime},x^{\prime})\\
\phi_{\mathbf{v}}(t,x)  & =\gamma\left[  \phi(t^{\prime},x^{\prime
})+\mathbf{v\cdot A}(t^{\prime},x^{\prime})\right] \\
\mathbf{A}_{\mathbf{v}}(t,x)  & =\gamma\left[  \mathbf{A}(t^{\prime}%
,x^{\prime})+\phi(t^{\prime},x^{\prime})\mathbf{v}\right]
\end{align*}
thus similarly to the case of eq. (\ref{KG}), in order to produce solitary
waves, it is sufficient to find stationary solutions of eq. (\ref{e1+}),
(\ref{e2+}), (\ref{e3+}) and to make a Lorentz transform.

In particular we shall look for solutions of (\ref{e1+}), (\ref{e2+}),
(\ref{e3+}) of the type
\begin{equation}
\psi(t,x)=u(x)e^{-i\omega t}\text{, }u\text{ }\geq0,\;\omega\text{ real,
}\mathbf{A=}0\mathbf{,}\text{ }\phi=\phi(x)\label{tipo}%
\end{equation}

We shall assume that $W$ is a $C^{2}$ function satisfying the following assumptions:

$W_{1})$ $W\geq0,$ $W(0)=W^{\prime}(0)=0$

$W_{2})W^{\prime\prime}(0)=m_{0}^{2}>0$

$W_{3})$ There exist $m_{1}$, $c>0$ with $m_{1}<m_{0}$ s.t.
\[
W(s)\leq\frac{1}{2}m_{1}^{2}s^{2}+c\text{ for all }s\in\mathbb{R}%
\]

$W_{4})$ for all $s\in\mathbb{R}$%
\[
0\leq\frac{1}{2}W^{\prime}(s)s\leq W(s)
\]

$W_{5})$ There exist constants $c_{1},c_{2}>0$ and $p<4$ such that for all $s
$%
\[
\left|  W^{\prime\prime}(s)\right|  \leq c_{1}\left|  s\right|  ^{p}+c_{2}%
\]

\bigskip

The main result of this paper is the following theorem:

\begin{theorem}
\label{gen}Assume that $W$ satisfies $W_{1}),...,W_{5})$. Then there exists
$q_{\ast}>0$ such that for any $q<q_{\ast}$ equations (\ref{e1+}),
(\ref{e2+}), (\ref{e3+}) possess a (non trivial) finite energy solution of the
type (\ref{tipo}).
\end{theorem}

\bigskip

\section{Existence of solitary waves in Abelian Gauge Theories}

The remaining part of this paper is devoted to the proof of a theorem of which
theorem \ref{gen} is an immediate consequence.

\subsection{Statement of the main theorem}

In this section we introduce some technical preliminaries and state the main theorem.

We look for solutions of (\ref{e1+}), (\ref{e2+}), (\ref{e3+}) of type
(\ref{tipo}). With this ansatz, equations (\ref{e2}) and (\ref{e4}) are
identically satisfied, while (\ref{e1}) and (\ref{e3}) become
\begin{equation}
-\Delta u-\left(  q\phi-\omega\right)  ^{2}u+W^{\prime}(u)=0\label{a11}%
\end{equation}
\begin{equation}
\Delta\phi=q\left(  q\phi-\omega\right)  u^{2}\;\label{a222}%
\end{equation}

We shall set
\[
\Phi=\frac{\phi}{\omega}%
\]
Equations (\ref{a11}), (\ref{a222}) become
\begin{equation}
-\Delta u-\omega^{2}\left(  q\Phi-1\right)  ^{2}\,u+W^{\prime}(u)=0\label{a12}%
\end{equation}
\begin{equation}
-\Delta\Phi+q^{2}u^{2}\Phi=qu^{2}\label{a23}%
\end{equation}

Let $H^{1}=$ $H^{1}(\mathbb{R}^{3})$ denote the usual Sobolev space with norm
\[
\left\Vert u\right\Vert _{H^{1}}=\left(  \int\left(  \left\vert \nabla
u\right\vert ^{2}+u^{2}\right)  dx\right)  ^{\frac{1}{2}}%
\]
and $D$ denote the completion of $C_{0}^{\infty}(\mathbb{R}^{3})$ with respect
to the inner product
\begin{equation}
\left(  v\mid w\right)  _{D}=\int\left(  \nabla v\mid\nabla w\right)  dx
\end{equation}
The request that $\psi(t,x)=u(x)e^{-i\omega t}$ and $\phi$ (or $\Phi$) possess
finite energy is satisfied if and only if $u\in H^{1}$ and $\phi$ (or $\Phi
$)$\in D.$

Now we can state the following theorem

\begin{theorem}
\textbf{(Main theorem)}\ \label{main}Assume that $W$ satisfies $W_{1}%
),...,W_{5}).$ Then there exists $q_{\ast}>0$ such that for any $q<q_{\ast}$
there exist $\omega\neq0$, $\omega^{2}<m_{0}^{2},$ and non trivial solutions
$u\in H^{1},$ $\Phi\in D$ of (\ref{a12}) and (\ref{a23}).
\end{theorem}

\subsection{The variational framework}

Consider the functional
\begin{equation}
F_{\omega}(u,\Phi)=J(u)-\omega^{2}\mathcal{A(}u,\Phi)\label{functional}%
\end{equation}
where
\begin{equation}
J(u)=\frac{1}{2}\int\left\vert \nabla u\right\vert ^{2}dx+\int
W(u)dx\label{new1}%
\end{equation}

and
\begin{equation}
\mathcal{A(}u,\Phi)=\frac{1}{2}\int\left\vert \nabla\Phi\right\vert
^{2}dx+\frac{1}{2}\int u^{2}(1-q\Phi)^{2}dx\label{new2}%
\end{equation}

Standard arguments show that $F_{\omega}$ is $C^{1}$ on $H^{1}(\mathbb{R}%
^{3})\times D$ and its critical points $(u,\Phi)$ are weak solutions of
(\ref{a12}) and (\ref{a23}); so equations (\ref{a12}) and (\ref{a23}) can be
written as follows
\begin{equation}
\frac{\partial F_{\omega}}{\partial u}(u,\Phi)=0\label{lak}%
\end{equation}

\begin{equation}
\frac{\partial F_{\omega}}{\partial\Phi}(u,\Phi)=0\text{ or }\frac
{\partial\mathcal{A}(u,\Phi)}{\partial\Phi}=0\label{lek}%
\end{equation}

where $\frac{\partial F_{\omega}}{\partial u}(u,\Phi),$ $\frac{\partial
F_{\omega}}{\partial\Phi}(u,\Phi),$ $\frac{\partial\mathcal{A}(u,\Phi
)}{\partial\Phi}$ denote the partial derivatives of $F_{\omega}$ and
$\mathcal{A}$ at $(u,\Phi)\in H^{1}(\mathbb{R}^{3})\times D,$ namely for any
$v\in H^{1}(\mathbb{R}^{3})$ and $w\in D$%
\begin{equation}
\frac{\partial F_{\omega}}{\partial u}(u,\Phi)\left[  v\right]  =\int\left(
(\nabla u\mid\nabla v)+W^{\prime}(u)v\right)  dx-\omega^{2}\int(1-q\Phi
)^{2}\,uvdx\label{tre}%
\end{equation}

\begin{align}
\frac{\partial F_{\omega}}{\partial\Phi}(u,\Phi)\left[  w\right]   &
=-\omega^{2}\frac{\partial\mathcal{A}(u,\Phi)}{\partial\Phi}(u,\Phi)\left[
w\right]  =\nonumber\\
& =-\omega^{2}\int\left(  (\nabla\Phi\mid\nabla w)+(q\Phi-1)\,u^{2}w\right)
dx\label{quattro}%
\end{align}

Now we reduce the study of (\ref{functional}) to the study of a functional of
the only variable $u$. Following \cite{bf}, it can be shown that, for any
$u\in H^{1}(\mathbb{R}^{3}),$ there exists a unique solution $\Phi\in D$ of
(\ref{a23}). So we can define the map%

\begin{equation}
u\in H^{1}(\mathbb{R}^{3})\rightarrow\Phi(u)\in D\label{map}%
\end{equation}

\noindent where $\Phi(u)$ is the unique solution of (\ref{a23}).

Now we set%

\begin{equation}
I_{\omega}(u)=F_{\omega}(u,\Phi(u))=J(u)-\omega^{2}\mathcal{A(}u,\Phi
(u)),\text{ }u\in H^{1}(\mathbb{R}^{3})\label{definition}%
\end{equation}
where $\Phi(u)$ is defined in (\ref{map}). The functional $I_{\omega}$ is
$C^{1}.$ Moreover, by (\ref{a23}), we easily get
\begin{equation}
\int\left|  \nabla\Phi(u)\right|  ^{2}dx=\int\left(  qu^{2}\Phi(u)-q^{2}%
u^{2}\Phi(u)^{2}\right)  dx\label{per}%
\end{equation}

Inserting (\ref{per}) in (\ref{new2}), we get%

\begin{equation}
\mathcal{A(}u,\Phi(u))=\frac{1}{2}\int\,u^{2}(1-q\Phi(u))dx\label{imp}%
\end{equation}

So by (\ref{definition}) and (\ref{imp}) we can write
\begin{equation}
I_{\omega}(u)=F_{\omega}(u,\Phi(u))=J(u)-\frac{\omega^{2}}{2}\int
\,u^{2}(1-q\Phi(u))dx\label{definition2}%
\end{equation}

The following proposition holds

\begin{proposition}
\label{cinque}Let $(u,\Phi)$ $\in H^{1}(\mathbb{R}^{3})\times D.$ Then the
following statements are equivalent:

a) $(u,\Phi)$ is a critical point of $F_{\omega}$

b) $u$ is a critical point of $I_{\omega}$ and $\Phi=\Phi(u)$ solves
(\ref{lek})
\end{proposition}%

\proof
Clearly we have
\[
\text{b)}\Leftrightarrow(\frac{\partial F_{\omega}}{\partial u}(u,\Phi
)+\frac{\partial F_{\omega}}{\partial\Phi}(u,\Phi)\Phi^{\prime}(u)\text{
}=0\text{ and }\Phi=\Phi(u)\text{ solves }(\ref{lek}))\Leftrightarrow
\]
\[
\Leftrightarrow(\frac{\partial F_{\omega}}{\partial u}(u,\Phi)=0\text{ and
}\frac{\partial F_{\omega}}{\partial\Phi}(u,\Phi)=0)\Leftrightarrow\text{a)}%
\]
%

\endproof

Then we are reduced to find the critical points of
\begin{equation}
I_{_{\omega}}(u)=J(u)-\frac{\omega^{2}}{2}\int\,u^{2}(1-q\Phi(u))dx,\text{
}u\in H^{1}(\mathbb{R}^{3})\label{princ}%
\end{equation}

The functional $I_{_{\omega}}$ presents a lack of compactness due to its
invariance under the group transformations $u(x)\rightarrow u(x+a)$
($a\in\mathbb{R}^{3}).$ To overcome this difficulty we restrict ourselves to
radial functions $u=u(r),$ $r=\left|  x\right|  .$ More precisely we shall
consider the functional $I_{_{\omega}}$ on the subspace of the radially
symmetric functions
\begin{equation}
H_{r}^{1}=\left\{  u\in H^{1}(\mathbb{R}^{3}):u=u(r),r=\left|  x\right|
\right\} \label{rad}%
\end{equation}

We recall (see \cite{strauss} or \cite{BL81}) that, for $6>p>2,$ $H_{r}^{1}$
is compactly embedded into $L_{r}^{p}$, where $L_{r}^{p}=\left\{  u\in
L^{p}(\mathbb{R}^{3}):u\text{ radially symmetric}\right\}  $.

Now $H_{r}^{1}$ is a natural constraint for $I_{\omega},$ namely the following
lemma holds

\begin{lemma}
\label{simm}Any critical point $u\in$ $H_{r}^{1}$ of $I_{_{\omega}}\mid
_{H_{r}^{1}}$ is also a critical point of $I_{\omega}.$
\end{lemma}%

\proof
Consider the $O(3)$ group action $T_{g}$ on $H^{1}(\mathbb{R}^{3})$ defined
by
\[
\text{for }g\in O(3),\text{ }u\text{ }\in H^{1}(\mathbb{R}^{3}):\text{ }%
T_{g}u(x)=u(g(x))
\]
Clearly $H_{r}^{1}$ is the set of the fixed points for this action, namely
\[
H_{r}^{1}=\left\{  u\in H^{1}(\mathbb{R}^{3})\mid T_{g}u=u\text{ for all }g\in
O(3)\right\}
\]
Then the conclusion can be achieved by usual arguments (see \cite{strauss}),
if we show that $I_{_{\omega}}$ is invariant under the $T_{g}$ action, namely
if
\begin{equation}
\text{for any }u\in H^{1}(\mathbb{R}^{3})\text{ and }g\in O(3)\text{ we have
}I_{_{\omega}}(T_{g}u)=I_{_{\omega}}(u)\label{inv}%
\end{equation}
Now, for $u\in H^{1}(\mathbb{R}^{3}),$ $\Phi(u)$ solves the equation
\[
-\Delta\Phi+q^{2}u^{2}\Phi=qu^{2}%
\]
then, if $g\in O(3),$ we have
\[
T_{g}(-\Delta\Phi(u)+q^{2}u^{2}\Phi(u))=qT_{g}(u^{2})
\]
\[
-\Delta(T_{g}\Phi(u))+q^{2}(T_{g}u)^{2}(T_{g}\Phi(u)))=q(T_{g}u)^{2}%
\]
This equality and the definition of $\Phi$ imply that
\begin{equation}
T_{g}\Phi(u)=\Phi(T_{g}u)\label{comm}%
\end{equation}

Therefore, using (\ref{comm}), we easily deduce (\ref{inv}).%

\endproof

We want to find non trivial $u,$ $\Phi,$ $\omega,$ such that (\ref{a12}) and
(\ref{a23}) are solved. Then, by lemma \ref{simm} and proposition
\ref{cinque}, we have to find for some $\omega\neq0$ a non trivial critical
point $u\in H_{r}^{1}$ of $I_{_{\omega}}\mid_{H_{r}^{1}}.$ In this case $u, $
$\Phi(u)$ will be non trivial weak solutions of (\ref{a12}) and (\ref{a23}).
To do this we set, for $\sigma>0$
\[
V_{\sigma}=\left\{  u\in H_{r}^{1}:\frac{1}{2}\int\,u^{2}(1-q\Phi
(u))dx=\sigma^{2}\right\}
\]
We shall look for critical points (in particular minimizers) of the functional
$J$ defined in (\ref{new1}) on $V_{\sigma}$ and $\omega^{2}$ will be obtained
as Lagrange multiplier.

Then in order to prove theorem \ref{main}, it will be enough to prove the
following proposition

\begin{proposition}
\label{min}Assume that $W$ satisfies the assumptions $W_{1}),...,W_{5}).$ Then
there exists $q_{\ast}>0$ such that for any $q<q_{\ast}$ there exists
$\sigma>0$ s.t. the functional $J$ has a minimizer on the manifold
\[
V_{\sigma}=\left\{  u\in H_{r}^{1}:\frac{1}{2}\int\,u^{2}(1-q\Phi
(u))dx=\sigma^{2}\right\}
\]
Moreover for the corresponding Lagrange multiplier $\omega^{2}$ we have
\[
0<\omega^{2}<m_{0}^{2}%
\]

\end{proposition}

\subsection{Preliminary lemmas}

In this section we shall prove some technical lemmas.

\begin{lemma}
\label{uno} The functional
\[
\Lambda:H^{1}\rightarrow\mathbb{R}^{3},\text{ }\Lambda(u)=\frac{1}{2}%
\int\,u^{2}(1-q\Phi(u))dx
\]
is $C^{1}$ and for any $u\in H^{1}$%
\begin{equation}
\Lambda^{\prime}(u)=u(1-q\Phi(u))\label{sec}%
\end{equation}
\bigskip
\end{lemma}%

\proof
Standard arguments show that $\Lambda$ is $C^{1}.$ Now we prove (\ref{sec}).
By (\ref{imp}) and (\ref{new2}) we have
\[
\Lambda(u)=\mathcal{A(}u,\Phi(u))=\frac{1}{2}\int\left\vert \nabla
\Phi(u)\right\vert ^{2}dx+\frac{1}{2}\int u^{2}(1-q\Phi(u))^{2}dx
\]
Then
\begin{equation}
\Lambda^{\prime}(u)=\frac{\partial\mathcal{A}}{\partial u}\mathcal{(}%
u,\Phi(u))+\frac{\partial\mathcal{A}}{\partial\Phi}\mathcal{(}u,\Phi
(u))\Phi^{\prime}(u)\label{give}%
\end{equation}

Since $\Phi(u)$ solves (\ref{a23}), we have
\[
\frac{\partial\mathcal{A}}{\partial\Phi}\mathcal{(}u,\Phi(u))=0
\]
Then (\ref{give}) gives
\[
\Lambda^{\prime}(u)=\frac{\partial\mathcal{A}}{\partial u}\mathcal{(}%
u,\Phi(u))=u(1-q\Phi(u))
\]
\endproof

The following lemma holds

\begin{lemma}
\label{massimo}Let $u\in H^{1}$ and $\Phi(u)\in D$ be the solution of
(\ref{a23}). Then
\begin{equation}
0\leq\Phi(u)\leq\frac{1}{q}%
\end{equation}

\end{lemma}%

\proof
Arguing by contradiction, we assume that there exists an open subset
$\Omega\subset\mathbb{R}^{3}$ such that
\begin{equation}
\Phi(u)>\frac{1}{q}\text{ in }\Omega\label{as}%
\end{equation}

\noindent then, since $\Phi(u)$ solves (\ref{a23}), we have
\[
-\Delta(\Phi(u)-\frac{1}{q})+q^{2}u^{2}(\Phi(u)-\frac{1}{q})=-\Delta
\Phi(u)+q^{2}u^{2}\Phi(u)-qu^{2}=0
\]
So \ $v=\Phi(u)-\frac{1}{q}$ satisfies
\[
-\Delta v+q^{2}u^{2}v=0\text{ in }\Omega,\text{ }v=0\text{ on }\partial\Omega
\]
then $v=0,$ contradicting (\ref{as}). An analogous argument shows that
$\Phi(u)\geq0.$%

\endproof

\begin{lemma}
\label{mani} For any $\sigma>0,$ the set
\[
V_{\sigma}=\left\{  u\in H_{r}^{1}:\frac{1}{2}\int\,u^{2}(1-q\Phi
(u))dx=\sigma^{2}\right\}
\]
is not empty and it is a one codimensional manifold.
\end{lemma}%

\proof
Let $\sigma>0$ and first prove that
\[
V_{\sigma}=\left\{  u\in H_{r}^{1}:\frac{1}{2}\int\,u^{2}(1-q\Phi
(u))dx=\sigma^{2}\right\}  \neq\emptyset.
\]
Fix $u\in H_{r}^{1}$ , $u\neq0$ a.e. in $\mathbb{R}^{3}$ and set for
$\lambda>0$%
\[
u_{\lambda}(x)=\lambda u(\lambda x),\text{ }\Phi_{\lambda}(u)(x)=\text{ }%
\Phi(u)(\lambda x)
\]
We have
\begin{equation}
\Phi_{\lambda}(u)(x)=\text{ }\Phi(u_{\lambda})(x)\label{prin}%
\end{equation}
In fact
\begin{multline}
-\Delta\Phi_{\lambda}(u)(x)+q^{2}u_{\lambda}^{2}(x)\Phi_{\lambda}(u)(x)=\\
=\lambda^{2}\left(  -\Delta\Phi(u)(\lambda x)+q^{2}u^{2}(\lambda
x)\Phi(u)(\lambda x)\right) \label{compl}%
\end{multline}

and, since $\Phi(u)$ satisfies (\ref{a23}), we have
\begin{equation}
-\Delta\Phi(u)(\lambda x)+q^{2}u^{2}(\lambda x)\Phi(u)(\lambda x)=qu^{2}%
(\lambda x)\label{coplex}%
\end{equation}
>From (\ref{compl}) and (\ref{coplex}) we get
\[
-\Delta\Phi_{\lambda}(u)(x)+q^{2}u_{\lambda}^{2}(x)\Phi_{\lambda
}(u)(x)=\lambda^{2}qu^{2}(\lambda x)=qu_{\lambda}^{2}(x)
\]
and this implies (\ref{prin}).

Now set
\[
\sigma_{\lambda}^{2}=\frac{1}{2}\int\,u_{\lambda}^{2}(1-q\Phi(u_{\lambda}))dx
\]
and first show that
\begin{equation}
\sigma_{1}^{2}=\frac{1}{2}\int\,u^{2}(1-q\Phi(u))dx>0\label{poss}%
\end{equation}

in fact, arguing by contradiction, assume that
\[
\int u^{2}(1-q\Phi(u))dx=0
\]
then, by lemma \ref{massimo} and since $u\neq0$ a.e. in $\mathbb{R}^{3},$ we
have
\[
\Phi(u)=\frac{1}{q}%
\]
contradicting the fact that $\Phi(u)\in D\subset L^{6}.$

Now we have
\[
\sigma_{\lambda}^{2}=\frac{1}{2}\int\,u_{\lambda}^{2}(1-q\Phi(u_{\lambda
}))dx=\text{(by (\ref{prin}}))=\frac{1}{2}\int\,u_{\lambda}^{2}(1-q\Phi
_{\lambda}(u))dx=
\]%
\[
=\frac{1}{2\lambda}\int\,u^{2}(1-q\Phi(u))dx=\frac{\sigma_{1}^{2}}{\lambda}%
\]

We conclude that for any $\sigma>0$ we can chose $\lambda=\frac{\sigma_{1}%
^{2}}{\sigma^{2}}$ so that $\sigma_{\lambda}^{2}=\sigma^{2}.$ This means that
$u_{\lambda}\in V_{\sigma}.$

$V_{\sigma}$ is a one codimensional manifold. In fact, by lemma \ref{uno},
$\Lambda$ is $C^{1}$ and we have
\[
\Lambda^{\prime}(u)=u(1-q\Phi(u))
\]
Then, since $\Phi(u)\leq\frac{1}{q},$ we have
\begin{align*}
\left(  \Lambda^{\prime}(u)=0\right)   &  \Longleftrightarrow\left(
u(1-q\Phi(u))=0\right)  \Rightarrow\left(  u^{2}(1-q\Phi(u))=0\right)
\Longrightarrow\\
&  \Longrightarrow(\int\left(  u^{2}(1-q\Phi(u)))dx=0\right)  \Rightarrow
\left(  u\notin V_{\sigma}\right)
\end{align*}
So for any $u\in V_{\sigma}$ we have
\[
\Lambda^{\prime}(u)\neq0
\]%
\endproof

\begin{lemma}
\label{quasi}Assume that $W$ satisfies $W_{3}).$ Then there exists $q_{\ast
}>0$ such that for $q<q_{\ast}$ there exists $\bar{u}\in H_{r}^{1},\bar{u}%
\neq0$ s.t.
\[
\text{ }\frac{J(\bar{u})}{\frac{1}{2}\int\bar{u}^{2}(1-q\Phi(\bar{u}%
))dx}<m_{0}^{2}%
\]

\end{lemma}%

\proof
Let us first consider the case $q=0$ and prove that there exists $\bar{u}\in
H_{r}^{1},\bar{u}\neq0$ s.t.
\begin{equation}
\text{ }\frac{J(\bar{u})}{\frac{1}{2}\int\bar{u}^{2}dx}<m_{0}^{2}\label{ee}%
\end{equation}

To this end let $\ u$ be a non trivial, smooth and radial function with
compact support $K$ and set for $\lambda>0$%
\[
u_{\lambda}=\lambda u(\frac{x}{\lambda}),\text{ }\sigma_{\lambda}^{2}=\frac
{1}{2}\int u_{\lambda}^{2}dx,\text{ }K_{\lambda}=\lambda K
\]
By assumption $W_{3})$ \ we have%

\begin{align*}
\frac{J(u_{\lambda})}{\sigma_{\lambda}^{2}}  & =\frac{\frac{1}{2}%
\int_{K_{\lambda}}\left\vert \nabla u_{\lambda}\right\vert ^{2}dx+\int
_{K_{\lambda}}W(u_{\lambda})dx}{\sigma_{\lambda}^{2}}\leq\\
& =\frac{\frac{1}{2}\int_{K_{\lambda}}\left\vert \nabla u_{\lambda}\right\vert
^{2}dx+\int_{K_{\lambda}}\left(  \frac{m_{1}^{2}}{2}u_{\lambda}^{2}+c\right)
dx}{\sigma_{\lambda}^{2}}=\\
& =\frac{\frac{1}{2}\int_{K_{\lambda}}\left\vert \nabla u_{\lambda}\right\vert
^{2}dx}{\sigma_{\lambda}^{2}}+\frac{\frac{m_{1}^{2}}{2}\int_{K_{\lambda}%
}u_{\lambda}^{2}dx}{\sigma_{\lambda}^{2}}+c\frac{meas(K_{\lambda})}%
{\sigma_{\lambda}^{2}}%
\end{align*}

Then by easy computations we get
\begin{equation}
\frac{J(u_{\lambda})}{\sigma_{\lambda}^{2}}\leq c_{1}\lambda^{-2}+m_{1}%
^{2}\label{be}%
\end{equation}
where $c_{1}$ is a positive constant depending only on the fixed map $u.$
Then, since $m_{1}^{2}<m_{0}^{2},$ by (\ref{be}) we get that (\ref{ee}) is
verified if we take
\[
\bar{u}=u_{\lambda}\text{, for }\lambda\text{ large}%
\]
Now we shall write $\Phi_{q}(\bar{u})$ instead of $\Phi(\bar{u})$ in order to
emphasize the dependence on $q.$ Clearly the conclusion will easily follow
from (\ref{ee}), if we show that
\begin{equation}
\frac{1}{2}\int\bar{u}^{2}(1-q\Phi_{q}(\bar{u}))dx\rightarrow\frac{1}{2}%
\int\bar{u}^{2}dx\text{ for }q\rightarrow0\label{prove}%
\end{equation}

So it remains to prove (\ref{prove}). The map $\Phi_{q}(\bar{u})$ satisfies
equation (\ref{a23}), then
\[
-\Delta\Phi_{q}(\bar{u})+q^{2}u^{2}\Phi_{q}(\bar{u})=q\bar{u}^{2}%
\]
So, multiplying by $\Phi_{q}(\bar{u})$ and integrating, we have
\[
\left\Vert \Phi_{q}(\bar{u})\right\Vert _{D}^{2}+q^{2}\int\bar{u}^{2}\Phi
_{q}(\bar{u})^{2}dx=q\int\bar{u}^{2}\Phi_{q}(\bar{u})dx\leq
\]%
\begin{equation}
\leq q\left\Vert \bar{u}\right\Vert _{L^{\frac{12}{5}}}^{2}\left\Vert \Phi
_{q}(\bar{u})\right\Vert _{L^{6}}\label{again}%
\end{equation}

and then
\[
\frac{\left\|  \Phi_{q}(\bar{u})\right\|  _{D}^{2}}{\left\|  \Phi_{q}(\bar
{u})\right\|  _{L^{6}}}\leq q\left\|  \bar{u}\right\|  _{L^{\frac{12}{5}}}^{2}%
\]
Then, since $D$ is continuously embedded into $L^{6}$, we easily get
\begin{equation}
\left\|  \Phi_{q}(\bar{u})\right\|  _{D}\leq c_{2}q\left\|  \bar{u}\right\|
_{L^{\frac{12}{5}}}^{2}\label{bi}%
\end{equation}

where $c_{2}$ ic a positive constant. Then, by (\ref{bi}) and by using again
(\ref{again}), we get
\[
q\int\bar{u}^{2}\Phi_{q}(\bar{u})dx\leq q\left\Vert \bar{u}\right\Vert
_{L^{\frac{12}{5}}}^{2}\left\Vert \Phi_{q}(\bar{u})\right\Vert _{L^{6}}\leq
c_{2}q^{2}\left\Vert \bar{u}\right\Vert _{L^{\frac{12}{5}}}^{4}%
\]
So we have
\[
q\int\bar{u}^{2}\Phi_{q}(\bar{u})dx\rightarrow0\text{ for }q\rightarrow0
\]
and (\ref{prove}) is proved.%

\endproof

\subsection{Proof of the main theorem}

Now we are ready to prove proposition \ref{min}.

Take $q<q_{\ast}$ and $\bar{u}$ as in lemma \ref{quasi} and set
\[
\sigma^{2}=\frac{1}{2}\int\bar{u}^{2}(1-q\Phi(\bar{u}))dx
\]
Consider
\[
V_{\sigma}=\left\{  u\in H_{r}^{1}:\frac{1}{2}\int\,u^{2}(1-q\Phi
(u))dx=\sigma^{2}\right\}
\]
We shall prove that $J\mid_{V_{\sigma}}$ has a\ minimizer whose Lagrange
multiplier $\omega^{2}$ satisfies
\[
0<\omega^{2}<m_{0}^{2}%
\]
Let $\left\{  u_{n}\right\}  \subset V_{\sigma}$ be a minimizing sequence for
$J\mid_{V_{\sigma}},$ i.e. $\left\{  u_{n}\right\}  \subset H_{r}^{1}$ s.t.
\[
\frac{1}{2}\int\,u_{n}^{2}(1-q\Phi(u_{n}))dx=\sigma^{2}%
\]
and
\[
J(u_{n})\rightarrow\inf J(V_{\sigma})
\]
By standard variational arguments we can assume that $\left\{  u_{n}\right\}
$ is also a criticizing sequence, i.e.
\begin{equation}
J^{\prime}\mid_{V_{\sigma}}(u_{n})\rightarrow0\label{ps}%
\end{equation}

Then there exists a sequence of real numbers $\left\{  \lambda_{n}\right\}  $
such that
\begin{equation}
J^{\prime}(u_{n})-\lambda_{n}\Lambda^{\prime}(u_{n})=\varepsilon
_{n}\rightarrow0\text{ in }H^{-1}\label{purre}%
\end{equation}

By \ (\ref{sec}) we can write
\begin{equation}
J^{\prime}(u_{n})-\lambda_{n}u_{n}(1-q\Phi_{n})=\varepsilon_{n}\rightarrow
0\text{ in }H^{-1}\label{bis}%
\end{equation}

where we have set
\[
\Phi_{n}=\Phi(u_{n})
\]
Now the study of the minimizing sequence $\left\{  u_{n}\right\}  $ is
essentially divided in two steps:

\begin{itemize}
\item First step: $\left\{  u_{n}\right\}  $ is bounded in $H^{1}.$

\item Second step: $\left\{  u_{n}\right\}  $ strongly converges (up to a
subsequence) in $H^{1}.$
\end{itemize}

Proof of the first step:

Clearly, since $W\geq0,$
\[
\left\{  \nabla u_{n}\right\}  \text{ is bounded in }L^{2}%
\]
So, in order to get the first step, we have only to show that $\left\{
u_{n}\right\}  $ is bounded in $L^{2}.$

We begin to show that
\begin{equation}
\left\{  \Phi_{n}\right\}  \text{is bounded in }D\label{first}%
\end{equation}

\bigskip Since the $\Phi_{n}^{\prime}$s solve (\ref{a23}) (with $u=u_{n}),$ we
get
\[
\int\left\vert \nabla\Phi_{n}\right\vert ^{2}dx=\int\left(  qu_{n}^{2}\Phi
_{n}-q^{2}u_{n}^{2}\Phi_{n}^{2}\right)  dx=
\]%
\[
=q\int u_{n}^{2}\Phi_{n}(1-q\Phi_{n})dx\leq(\text{by lemma \ref{massimo})}\leq
\]%
\[
\leq\int u_{n}^{2}(1-q\Phi_{n})dx=2\sigma^{2}%
\]
Then (\ref{first}) is verified.

The maps $\Phi_{n}$ are radially symmetric and the sequence $\left\{  \Phi
_{n}\right\}  $ is bounded in $D.$ Then, by virtue of a well known radial
lemma \cite{BL81}, we have
\[
\left|  \Phi_{n}(x)\right|  \leq c_{1}\left|  x\right|  ^{-\frac{1}{2}%
}\left\|  \Phi_{n}\right\|  _{D}\leq c_{2}\left|  x\right|  ^{-\frac{1}{2}%
}\text{ for }\left|  x\right|  \geq1
\]
where $c_{1},c_{2}$ are positive constants independent on $n.$

Then there exists $R>0$ large enough so that for all $n\in\mathbf{N}$ and for
$\left\vert x\right\vert \geq R$
\[
1-q\Phi_{n}(x)\geq k>0\text{ }%
\]
So
\begin{equation}
2\sigma^{2}=\int_{{}}u_{n}^{2}(1-q\Phi_{n})dx\geq\int_{B_{R}^{c}}u_{n}%
^{2}(1-q\Phi_{n})dx\geq k\int_{B_{R}^{c}}u_{n}^{2}dx\label{conc}%
\end{equation}

where
\[
B_{R}=\left\{  x\in\mathbb{R}^{3}:\left|  x\right|  <R\right\}  ,\text{ }%
B_{R}^{c}=\mathbb{R}^{3}-B_{R}%
\]

\bigskip On the other hand, since $\left\{  u_{n}\right\}  $ is bounded in
$D\subset L^{6},$ we have that
\begin{equation}
\left\{  u_{n}\right\}  \text{is bounded in }L^{2}(B_{R})\label{conclu}%
\end{equation}

\bigskip Finally by (\ref{conc}) and (\ref{conclu}) we conclude that $\left\{
u_{n}\right\}  $is bounded in $L^{2}$

Now we prove the second step:

First we show that, up to subsequence,
\begin{equation}
\lambda_{n}\rightarrow\omega^{2}\text{ where \ \ \ \ \ \ \ \ }\omega^{2}%
<m_{0}^{2}\label{bou}%
\end{equation}

By (\ref{bis}) we get
\begin{equation}
\left\langle J^{\prime}(u_{n}),u_{n}\right\rangle -\lambda_{n}\int u_{n}%
^{2}(1-q\Phi_{n})dx=\left\langle J^{\prime}(u_{n}),u_{n}\right\rangle
-2\lambda_{n}\sigma^{2}=\delta_{n}\label{cosi}%
\end{equation}

where $\delta_{n}=\left\langle \varepsilon_{n},u_{n}\right\rangle $ and
$\left\langle \text{ \ \ \ }\right\rangle $ denotes the pairing between $H^{1}
$ and its dual $H^{-1}$. Since $\left\{  u_{n}\right\}  $ is bounded in
$H^{1}$ and $\varepsilon_{n}\rightarrow0$ in $H^{-1},$ we have
\begin{equation}
\delta_{n}\rightarrow0\label{coso}%
\end{equation}
Then
\begin{align}
\lambda_{n}  & =\frac{1}{2\sigma^{2}}\left(  \left\langle J^{\prime}%
(u_{n}),u_{n}\right\rangle -\delta_{n}\right)  =\nonumber\\
& =\frac{1}{\sigma^{2}}\left(  \frac{1}{2}\int\left(  \left\vert \nabla
u_{n}\right\vert ^{2}+W^{\prime}(u_{n})u_{n}\right)  dx\right)  -\frac
{\delta_{n}}{2\sigma^{2}}\leq\nonumber\\
(\text{by }W_{4})  & \leq\frac{J(u_{n})}{\sigma^{2}}-\frac{\delta_{n}}%
{2\sigma^{2}}.\label{pari}%
\end{align}

On the other hand
\begin{equation}
\left\langle J^{\prime}(u_{n}),u_{n}\right\rangle =\int\left(  \left\vert
\nabla u_{n}\right\vert ^{2}+W^{\prime}(u_{n})u_{n}\right)  dx\geq
0\label{ouno}%
\end{equation}

Then by (\ref{pari}) and (\ref{ouno}) we easily get
\begin{equation}
-\frac{\delta_{n}}{2\sigma^{2}}\leq\lambda_{n}\leq\frac{J(u_{n})}{\sigma^{2}%
}-\frac{\delta_{n}}{2\sigma^{2}}\text{ }\label{cosa}%
\end{equation}
Then, up to a subsequence, we have
\begin{equation}
\lambda_{n}\rightarrow\omega^{2}\label{dis}%
\end{equation}
where
\begin{equation}
0\leq\omega^{2}\leq\inf\left\{  \frac{J(u)}{\sigma^{2}}:u\in V_{\sigma
}\right\} \label{impo}%
\end{equation}
Since $\bar{u}$ and $q$ are chosen as in lemma \ref{quasi}, we have
\begin{equation}
\frac{J(\bar{u})}{\sigma^{2}}<m_{0}^{2}\label{impor}%
\end{equation}
From (\ref{dis}), (\ref{impo}) and (\ref{impor}) we clearly get (\ref{bou}).

Now we show that $\left\{  u_{n}\right\}  $ is a Cauchy sequence in $H^{1}.$
Set
\[
w_{nm}=u_{n}-u_{m}.
\]
Then writing (\ref{bis}) for $u_{n}$and $u_{m}$ and subtracting we get
\begin{equation}
-\Delta w_{nm}+W^{\prime}(u_{n})-W^{\prime}(u_{m})-\lambda_{n}u_{n}%
+\lambda_{m}u_{m}+q\lambda_{n}u_{n}\Phi_{n}-q\lambda_{m}u_{m}\Phi
_{m}=\varepsilon_{nm}\label{bene}%
\end{equation}

where
\[
\varepsilon_{nm}=\varepsilon_{n}-\varepsilon_{m}%
\]
If we add and subtract in (\ref{bene}) the terms $\lambda_{n}u_{m},$
$q\lambda_{n}u_{m}\Phi_{m}$ and multiply, by the duality map, both the sides
by $w_{nm}$, we get
\begin{multline}
\int\left(  \left\vert \nabla w_{nm}\right\vert ^{2}+\left(  W^{\prime}%
(u_{n})-W^{\prime}(u_{m})\right)  w_{nm}-\lambda_{n}w_{nm}^{2}+(\lambda
_{m}-\lambda_{n})u_{m}w_{nm}\right)  dx+\\
+q(\lambda_{n}-\lambda_{m})\int u_{m}\Phi_{m}w_{nm}dx+q\lambda_{n}\int\left(
u_{n}\Phi_{n}-u_{m}\Phi_{m}\right)  w_{nm}dx=\\
=\left\langle \varepsilon_{nm},w_{nm}\right\rangle =\delta_{nm}\label{pre}%
\end{multline}

Since \bigskip$w_{nm}$ are bounded in $H^{1}$ and $\varepsilon_{nm}%
\rightarrow0$ in $H^{-1},$ we get
\begin{equation}
\delta_{nm}\rightarrow0\label{for}%
\end{equation}
Now
\begin{equation}
\int\left(  W^{\prime}(u_{n})-W^{\prime}(u_{m})\right)  w_{nm}dx=\int\left(
W^{\prime\prime}(\xi_{nm})\right)  w_{nm}^{2}dx\label{fer}%
\end{equation}

where
\begin{equation}
\xi_{nm}=tu_{n}+(1-t)u_{m}\text{, }0\leq t\leq1\label{fur}%
\end{equation}
Therefore (\ref{pre}) can be written as follows
\begin{equation}
\int\left(  \left|  \nabla w_{nm}\right|  ^{2}+\left(  W^{\prime\prime}%
(\xi_{nm})-\lambda_{n}\right)  w_{nm}^{2}\right)  dx+B_{nm}^{1}+B_{nm}%
^{2}+B_{nm}^{3}=\delta_{nm}\label{written}%
\end{equation}

where
\[
B_{nm}^{1}=\int(\lambda_{m}-\lambda_{n})u_{m}w_{nm}dx,\text{ }B_{nm}%
^{2}=q(\lambda_{n}-\lambda_{m})\int u_{m}\Phi_{m}w_{nm}dx
\]
and
\[
B_{nm}^{3}=q\lambda_{n}\int\left(  u_{n}\Phi_{n}-u_{m}\Phi_{m}\right)
w_{nm}dx
\]
Now we show that
\begin{equation}
B_{nm}^{1}\rightarrow0,B_{nm}^{2}\rightarrow0,B_{nm}^{3}\rightarrow
0\label{zero}%
\end{equation}
Since $\left\{  u_{n}\right\}  $ is bounded in $H^{1},$ the sequence $\left\{
\int u_{m}w_{nm}\right\}  $ is also bounded. Moreover $\lambda_{n}-\lambda
_{m}\rightarrow0,$ then
\[
B_{nm}^{1}\rightarrow0
\]
Analogously, in order to prove that \ $B_{nm}^{2}\rightarrow0,$ it will be
enough to show that
\begin{equation}
\left\{  \int u_{m}\Phi_{m}w_{nm}dx\right\}  \text{is bounded}\label{true}%
\end{equation}
This is true, in fact
\[
\left\vert \int u_{m}\Phi_{m}w_{nm}dx\right\vert \leq\left\Vert u_{m}%
\right\Vert _{L^{\frac{12}{5}}}\left\Vert \Phi_{m}\right\Vert _{L^{6}%
}\left\Vert w_{nm}\right\Vert _{L^{\frac{12}{5}}}%
\]
and we know that $\left\{  u_{m}\right\}  $ is bounded in $H^{1}($ and
therefore in $L^{\frac{12}{5}})$ and $\left\{  \Phi_{m}\right\}  $ is bounded
in $D$ (and therefore in $L^{6}).$

Finally in order to prove that $B_{nm}^{3}\rightarrow0,$ it will be enough to
show that
\begin{equation}
\int\left(  u_{n}\Phi_{n}-u_{m}\Phi_{m}\right)  w_{nm}dx\rightarrow
0\label{finas}%
\end{equation}
Clearly
\begin{equation}
\int\left(  u_{n}\Phi_{n}\right)  w_{nm}dx\leq\left\Vert w_{nm}\right\Vert
_{L^{3}}\left\Vert u_{n}\Phi_{n}\right\Vert _{L^{\frac{3}{2}}}\label{an}%
\end{equation}

Since $H_{r}^{1}$ is compactly embedded into $L^{3}$ (see \cite{BL81} or
\cite{strauss})$,$ and $w_{nm}=u_{n}-u_{m}$ weakly converges to $0$ in
$H_{r}^{1},$ we have
\begin{equation}
w_{nm}=u_{n}-u_{m}\rightarrow0\text{ in }L^{3}\label{on}%
\end{equation}
On the other hand
\[
\int\left(  u_{n}\Phi_{n}\right)  ^{\frac{3}{2}}dx\leq\left(  \int u_{n}%
^{2}dx\right)  ^{\frac{3}{4}}\left(  \int\Phi_{n}^{6}dx\right)  ^{\frac{1}{4}}%
\]
Then, since $\left\{  u_{m}\right\}  $ is bounded in $L^{2}$ and $\left\{
\Phi_{n}\right\}  $ is bounded in $L^{6},\left\Vert u_{n}\Phi_{n}\right\Vert
_{L^{\frac{3}{2}}}$ is bounded. So (\ref{finas}) easily follows from
(\ref{an}) and (\ref{on})..

We conclude that (\ref{zero}) holds and therefore (\ref{written}) implies
that
\begin{equation}
\int\left(  \left\vert \nabla w_{nm}\right\vert ^{2}+\left(  W^{\prime\prime
}(\xi_{nm})-\lambda_{n}\right)  w_{nm}^{2}\right)  dx=\eta_{nm}\label{mu}%
\end{equation}

where $\eta_{nm}\rightarrow0.$

Now, since $\lambda_{n}\rightarrow\omega^{2}<m_{0}^{2}$ (see (\ref{bou}) and
$W^{\prime\prime}(0)=m_{0}^{2},$ there exist $\theta,\delta>0$ such that
\begin{equation}
W^{\prime\prime}(s)-\lambda_{n}>\theta\text{ \ for }n\text{ large and
}\left\vert s\right\vert <\delta\label{ma}%
\end{equation}

By a well known radial lemma (\cite{BL81} or \cite{strauss}) and since the
sequence $\left\{  \xi_{nm}\right\}  $ is bounded in $H_{r}^{1},$ there exists
$M>0$ such that for all $n,m$
\begin{equation}
\left|  \xi_{nm}(x)\right|  <\delta\text{ for }\left|  x\right|  >M\text{
}\label{mi}%
\end{equation}
Then, if we set
\[
B_{M}=\left\{  x\in\mathbb{R}^{3}:\left|  x\right|  <M\right\}  ,\text{ }%
B_{M}^{c}=\mathbb{R}^{3}-B_{M},
\]
from (\ref{ma}) and (\ref{mi}) we easily get for $n$ large
\begin{equation}
\int_{B_{M}^{c}}\left(  W^{\prime\prime}(\xi_{nm})-\lambda_{n}\right)
w_{nm}^{2}dx\geq\int_{B_{M}^{c}}\theta w_{nm}^{2}dx\label{fio}%
\end{equation}

Now we show that
\begin{equation}
\int_{B_{M}}\left(  W^{\prime\prime}(\xi_{nm})-\lambda_{n})w_{nm}^{2}\right)
dx\rightarrow0\label{fai}%
\end{equation}
The Sobolev space $H^{1}(B_{M})$ is compactly embedded into $L^{q}(B_{M})$
($2\leq q<6),$ and the sequence $w_{nm}=u_{n}-u_{m}$ weakly converges to $0$
in $H^{1}(B_{M}).$ Then we have
\begin{equation}
\left\Vert w_{nm}\right\Vert _{L^{q}(B_{M})}\rightarrow0\text{ for }2\leq
q<6\label{ancor}%
\end{equation}
By assumption $W_{5})$ we easily get
\[
\int_{B_{M}}\left\vert W^{\prime\prime}(\xi_{nm})w_{nm}^{2}\right\vert dx\leq
c_{1}\int_{B_{M}}\left\vert \xi_{nm}\right\vert ^{p}w_{nm}^{2}dx+c_{2}%
\int_{B_{M}}w_{nm}^{2}dx\leq
\]%
\begin{equation}
c_{1}\left\Vert \xi_{nm}\right\Vert _{L^{6}(B_{M})}^{p}\left\Vert
w_{nm}\right\Vert _{L^{\bar{q}}(B_{M})}^{2}+c_{2}\left\Vert w_{nm}\right\Vert
_{L^{2}(B_{M})}^{2}\label{eric}%
\end{equation}

where
\[
\bar{q}=\frac{12}{6-p}<6
\]
Since the sequence $\left\{  \left\Vert \xi_{nm}\right\Vert _{L^{6}(B_{M}%
)}^{p}\right\}  $ is bounded, we deduce by (\ref{ancor}) and (\ref{eric}) that%
\[
\int_{B_{M}}\left\vert W^{\prime\prime}(\xi_{nm})w_{nm}^{2}\right\vert
dx\rightarrow0
\]
Then we clearly have (\ref{fai}).

From (\ref{ancor}), (\ref{fio}) and (\ref{fai}) we easily get
\begin{equation}
\int\left(  W^{\prime\prime}(\xi_{nm})-\lambda_{n})w_{nm}^{2}\right)
dx\geq\int\theta w_{nm}^{2}dx+\mu_{nm},\text{ with }\mu_{nm}\rightarrow
0\label{pil}%
\end{equation}

From (\ref{mu}) and (\ref{pil}) we get
\[
\int\left(  \left\vert \nabla w_{nm}\right\vert ^{2}+\theta w_{nm}^{2}\right)
dx\leq\gamma_{nm}%
\]
where $\gamma_{nm}=$\ $\eta_{nm}-$\ $\mu_{nm}\rightarrow0.$\ \ Then we
conclude that $\left\{  u_{n}\right\}  $ is a Cauchy sequence in $H^{1}$ and
therefore
\begin{equation}
u_{n}\rightarrow u\text{ strongly in }H^{1}\label{str}%
\end{equation}
Then, since $J$ and $\Lambda$ are $C^{1},$ we have
\[
J(u_{n})\rightarrow J(u)=\inf J(V_{\sigma})
\]
and
\begin{equation}
J^{\prime}(u_{n})\rightarrow J^{\prime}(u),\text{ }\Lambda^{\prime}%
(u_{n})\rightarrow\Lambda^{\prime}(u)\label{purr}%
\end{equation}

Now we prove that $u\in V_{\sigma},$ i.e.
\begin{equation}
\frac{1}{2}\int\,u^{2}(1-q\Phi(u))dx=\sigma^{2}\label{finni}%
\end{equation}

Since $\left\{  \Phi_{m}\right\}  $ is bounded in $D,$ we have (up to a
subsequence)
\begin{equation}
\Phi_{n}\rightharpoonup\Phi\text{ weakly in }D\label{we}%
\end{equation}

First observe that
\begin{equation}
\Phi=\Phi(u)\label{ho}%
\end{equation}

i.e. $\Phi$ solves the equation
\begin{equation}
-\Delta\Phi+q^{2}u^{2}\Phi=qu^{2}\label{eq}%
\end{equation}

In fact, since $\Phi_{n}=\Phi(u_{n}),$ we have
\begin{equation}
-\Delta\Phi_{n}+q^{2}u_{n}^{2}\Phi_{n}=qu_{n}^{2}\label{en}%
\end{equation}
Then, by (\ref{str}) and (\ref{we}), standard calculations permit to take the
limit in (\ref{en}) and to get (\ref{eq}). So (\ref{ho}) holds.

Since $u_{n}\in V_{\sigma},$ we have
\[
\frac{1}{2}\int\,u_{n}^{2}(1-q\Phi_{n})dx=\sigma^{2}\text{ for all }n
\]
Again, by (\ref{str}) and (\ref{we}), standard calculations permit to show
that
\[
\frac{1}{2}\int\,u_{n}^{2}(1-q\Phi_{n})dx\rightarrow\frac{1}{2}\int
\,u^{2}(1-q\Phi)dx=\frac{1}{2}\int\,u^{2}(1-q\Phi(u))dx
\]
So (\ref{finni}) holds and $u$ is a minimizer of $J$ on $V_{\sigma}$.

By (\ref{purr}), (\ref{purre}) and (\ref{bou}) $u$ solves the equation
\[
J^{\prime}(u)-\omega^{2}\Lambda^{\prime}(u)=0
\]
where $\omega^{2}<m_{0}^{2}.$ $.$

Finally we show that $\omega^{2}>0$. In fact, arguing by contradiction, assume
$\omega^{2}=0,$ then $u\in H^{1}$ is a non trivial solution of the equation
\[
J^{\prime}(u)=0
\]
i.e.
\[
-\Delta u+W^{\prime}(u)=0
\]
However, since $W\geq0,$ this equation has no non trivial solution in $H^{1} $
by the Derrick-Pohozaev identity (\ref{ventuno}).%

\endproof

$.$

\end{document}